\documentclass[11pt,twoside]{article}%
\usepackage{latexsym}
\usepackage{amsmath}
\usepackage{graphicx}%
\usepackage{amsfonts}%
\usepackage{amssymb}
\begin{document}

\title{Generating uniform random vectors in $\mathbf{Z}_{p}^{k}$: the general case}
\author{Claudio Asci\\Dipartimento di Matematica e Informatica\\{Universit\`{a} degli Studi di Trieste}}
\date{}
\maketitle

\begin{abstract}
This paper is about the rate of convergence of the Markov chain $\mathbf{X}%
_{n+1}=A\mathbf{X}_{n}+\mathbf{B}_{n}$ (mod $p$), where $A$ is an integer
matrix with nonzero eigenvalues and $\left\{  \mathbf{B}_{n}\right\}  _{n}$ is
a sequence of independent and identically distributed integer vectors, with
support not parallel to a proper subspace of $\mathbf{Q}^{k}$ invariant under
$A$. If $|\lambda_{i}|\not =1$ for all eigenvalues $\lambda_{i}$ of $A$, then
$n=O\left(  (\ln p)^{2}\right)  $ steps are sufficient and $n=O(\ln p)$ steps
are necessary to have $\mathbf{X}_{n}$ sampling from a nearly uniform
distribution. Conversely, if $A$ has the eigenvalues $\lambda_{i}$ that are
roots of positive integer numbers, $|\lambda_{1}|=1$ and $|\lambda_{i}|>1$ for
all $i\not =1$, then $O\left(  p^{2}\right)  $ steps are necessary and sufficient.

\end{abstract}

\noindent

\textbf{Running head.} Generating uniform random vectors.\newpage

\section{Introduction}

In this paper we generalize some results obtained in the paper \cite{Asci}
about Markov chains on $\mathbf{Z}^{k}$\ of the form%
\begin{equation}
\mathbf{X}_{n+1}=A\mathbf{X}_{n}+\mathbf{B}_{n}\ (\text{mod }p), \label{iter}%
\end{equation}
where $\mathbf{X}_{0}=\mathbf{x}_{0}\in\mathbf{Z}^{k}$, $A\in GL_{k}%
(\mathbf{Q})\cap M_{k}(\mathbf{Z})$, $p$ is an integer and $\left\{
\mathbf{B}_{n}\right\}  _{n}$ is a sequence of independent and identically
distributed integer vectors.

If $k=1$ and $\mathbf{B}_{n}\ $is a fixed integer $b$, for particular values
of $p$\ this recursion is used to produce pseudorandom numbers on computers
(see, for example, Knuth's book \cite{Knuth}).

In the paper \cite{Chung}, the constant term $b$ is chosen with a fixed
probability at each step and the authors study the following Markov chain:%

\[
X_{n+1}=aX_{n}+B_{n}\ (\text{mod }p),
\]
where $a$ is a positive integer. This randomness is introduced in order to
produce uniformly distributed random numbers on the set $\{0,1,...,p-1\}$. In
the cited paper, it is shown that, for $a=2$, $n=O(\ln p\ \ln\ln p)$ steps are
sufficient to sample $X_{n}$ from an almost uniform distribution. On the other
hand, if $a=1$ then $n=O\left(  p^{2}\right)  $ steps are necessary and
sufficient to achieve randomness. A further generalization is described in
\cite{Hildebrand-93}, where the integer $a$ is also allowed to vary.

The extension of the previous results to the higher-dimensional case is due to
Asci \cite{Asci} and next to Hildebrand and McCollum \cite{Hildebrand-07},
with the study of some particular cases of the recursion (\ref{iter}). In
\cite{Asci}, the distribution of $\mathbf{B}_{n}$ is the most general (the
support of the distribution cannot be parallel to any proper subspace of
$\mathbf{Q}^{k}$ invariant under $A$), but the matrix $A$ has only integer
eigenvalues. In \cite{Hildebrand-07}, $A$ is arbitrary, but only a specific
distribution for $\mathbf{B}_{n}$\ is considered.

The general case is studied in this paper, with the condition $\left\|
\mathbf{B}_{n}\right\|  _{\infty}\in L^{2}$ and some further conditions on
$p$. We find two different types of behaviour for the sequence (\ref{iter}),
depending on the size of the complex eigenvalues of $A$. If $|\lambda_{i}%
|\neq1$ for all eigenvalues $\lambda_{i}$, then $n=O\left(  (\ln
p)^{2}\right)  $ steps are sufficient and $n=O(\ln p)$ steps are necessary to
reach the uniform distribution (theorems 3.1 and 3.13). In particular, for a
matrix $A$ with eigenvalues $\lambda_{i}$ that are roots of positive integers
and $|\lambda_{i}|>1$, we can show that $n=O(\ln p\ln\ln p)$ steps are
sufficient (theorem 3.7). On the other hand, if the eigenvalue $\lambda_{i}$
are roots of positive integers, $|\lambda_{1}|=1$ and $|\lambda_{i}|>1$ for
all $i\not =1$, then $O\left(  p^{2}\right)  $ steps are necessary and
sufficient (theorems 3.12 and 3.14). These theorems agree with the
one-dimensional case studied in \cite{Chung} and with the results in
\cite{Asci} and \cite{Hildebrand-07}.

In Section 2, we recall some general results about random walk on groups and
the preliminary lemmas proved in \cite{Asci}. The main results of our work can
be found in Section 3.

\section{Preliminary results}

Consider the sequence (\ref{iter}) and observe that we can suppose
$\mathbf{X}_{n}\in\mathbf{Z}_{p}^{k}$.

Set $P_{n}(\mathbf{x})=p(\mathbf{X}_{n}=\mathbf{x})$,$\ \forall~\mathbf{x}%
\in\mathbf{Z}_{p}^{k}$, and $\mu(\mathbf{x})=p(\mathbf{B}_{n}=\mathbf{x}%
)$,$\ \forall~\mathbf{x}\in\mathbf{Z}^{k}$,$\ \forall~n\in\mathbf{N}$;
moreover, denote by $U$ the uniform distribution on $\mathbf{Z}_{p}^{k}$.
Define:
\[
V=\{\mathbf{x}\in\mathbf{Z}^{k}:\mathbf{x}=\mathbf{h}-\mathbf{k}%
,\quad\text{where}\ \mathbf{h},\mathbf{k}\in\text{supp }\mu\}.
\]
Indicate by $d$, where $d\leq k$, the degree of the minimum polynomial of $A$.
By definition:%
\[
\overset{d}{\underset{i=1}{\prod}}(A-\lambda_{i}I)=\overset{d}{\underset
{i=1}{\prod}}(^{t}A-\lambda_{i}I)=0\in M_{k}(\mathbf{Z}),\quad\lambda_{i}%
\in\{\lambda_{1},...,\lambda_{d}\},\ \forall~i=d+1,...,k,
\]
where $\lambda_{1},...,\lambda_{d},...,\lambda_{k}$ are the eigenvalues of
$A$. Finally, set:
\[
V^{d-1}=\{A^{m}\mathbf{x}:\mathbf{x}\in V\ ,m=0,1,...,d-1\}.
\]

In order to show that the distribution $P_{n}$ tends to the uniform
distribution $U$, as $n\rightarrow+\infty$, use the Fourier analysis (see
Diaconis' monograph \cite{Diaconis} and \cite{Helleloid}). Define the
variation distance between $P_{n}$ and $U$ in the following way:
\[
\Vert P_{n}-U\Vert={\frac{1}{2}}\sum_{\mathbf{\alpha}\in\mathbf{Z}_{p}^{k}%
}|P_{n}(\mathbf{\alpha})-U(\mathbf{\alpha})|.
\]

It is possible to prove that
\[
\Vert P_{n}-U\Vert={\frac{1}{2}}\sup_{f\in F}|E_{P_{n}}(f)-E_{U}%
(f)|=\max_{A\subset\mathbf{Z}_{p}^{k}}|P_{n}(A)-U(A)|,
\]
where $F\equiv\{f:\mathbf{Z}_{p}^{k}\longrightarrow\mathbf{C:}\Vert f\Vert
\leq1\mathbf{\}}$.

Henceforth, our purpose will be to find a bound for $\Vert P_{n}-U\Vert$ in
terms of $n$ and $p$. Observe that, if we indicate with $\left\{
\mathbf{Y}_{n}\right\}  _{n}$ the sequence defined by (\ref{iter}) and the
condition $\mathbf{X}_{0}=\mathbf{0}$, we have $\displaystyle\mathbf{X}%
_{n}=\varphi_{n}(\mathbf{Y}_{n})$, where the one to one function $\varphi
_{n}:\mathbf{Z}_{p}^{k}\longrightarrow\mathbf{Z}_{p}^{k}$ is defined by
$\displaystyle\varphi_{n}(\mathbf{x})=A^{n}\mathbf{x}_{0}+\mathbf{x}$.
Moreover:%
\[
\Vert P_{n}-U\Vert=\Vert(P_{n}\circ\varphi_{n})-U\Vert,
\]
then we can consider $\mathbf{X}_{0}=\mathbf{0}$.

Let $f:\mathbf{Z}^{k}\longrightarrow\mathbf{C}$; define the Fourier
transform$\ \widehat{{f}}:\mathbf{R}^{k}\longrightarrow C$ by:
\[
\widehat{{f}}(\mathbf{\alpha})=\sum_{\mathbf{h}\in\mathbf{Z}^{k}}\exp\left(
{\frac{2\pi i}{p}}\langle\mathbf{h},\mathbf{\alpha}\rangle\right)
f(\mathbf{h}).
\]

We have the following four results, whose proofs are similar to the proofs of
the lemmas 2.5, 3.1, 3.3 and 3.4 in \cite{Asci}: the only difference is that
$\mathbf{\alpha}$ ranges in $\mathbf{R}^{k}$ instead of in $\mathbf{Z}_{p}%
^{k}$. Lemma 2.1 is also proved in \cite{Diaconis}, in a more general case.

\textbf{Lemma 2.1} (Upper bound lemma)\textbf{.}
\[
\Vert P_{n}-U\Vert^{2}\leq{\frac{1}{4}}\sum_{\alpha\in\mathbf{Z}_{p}%
^{k}-\{\mathbf{0}\}}|\widehat{{P}}_{n}(\mathbf{\alpha})|^{2}.
\]
\medskip

\textbf{Lemma 2.2.}\ Suppose $\gcd(\det(A),p)=1$, $\mathbf{X}_{0}=\mathbf{0}$,
$\mathbf{\alpha}\in\mathbf{R}^{k}$; then:
\[
\widehat{{P}}_{n}(\mathbf{\alpha})=\prod_{j=0}^{n-1}\widehat{{\mu}}\left(
{}^{t}A^{j}\mathbf{\alpha}\right)  .\leqno1)
\]%
\[
|\widehat{{P}}_{n}(\mathbf{\alpha})|^{2}=\prod_{j=0}^{n-1}\left(
{\sum\limits_{\mathbf{h,i}\in\mathbf{Z}^{k}}{\mu(\mathbf{h})\mu(\mathbf{i}%
)\cos\left(  {\frac{{2\pi}}{p}}\left\langle \mathbf{h-i},{}^{t}A^{j}%
\mathbf{\alpha}\right\rangle \right)  }}\right)  \leqno2)
\]%
\[
\leq\prod_{j=0}^{n-1}\left(  {1-2\mu\left(  \mathbf{u}\right)  \mu\left(
\mathbf{v}\right)  +2\mu\left(  \mathbf{u}\right)  \mu\left(  \mathbf{v}%
\right)  \cos\left(  {\frac{{2\pi}}{p}\left\langle \mathbf{u-v},{}^{t}%
A^{j}\mathbf{\alpha}\right\rangle }\right)  }\right)  ,
\]
$\forall\,\mathbf{u,v}\in$ supp $\mu$.\medskip

\textbf{Lemma 2.3.}\ Let $\mathbf{\alpha}\in\mathbf{Z}_{p}^{k}-\{\mathbf{0\}}%
$; then:
\[
\Vert P_{n}-U\Vert\geq{\frac{1}{2}}\left|  \widehat{{P}}_{n}(\mathbf{\alpha
})\right|  .
\]
\medskip

\textbf{Lemma 2.4.}\ Suppose that the support of $\mu$ is not parallel to a
proper subspace of $\mathbf{Q}^{k}$ invariant under $A$. Then, there exists a
basis $\{\mathbf{y}_{1},...,\mathbf{y}_{k}\}\subset V^{d-1}$ of $\mathbf{Q}%
^{k}$. Furthermore, for all $p\in\mathbf{N}$ such that $\gcd(\det
(\mathbf{y}_{1},...,\mathbf{y}_{k}),p)=1$, for all $\mathbf{\alpha}%
\in\mathbf{R}^{k}-(p\mathbf{Z)}^{k}$, there exists $i\in\{1,...,k\}$ such that
$\langle\mathbf{y}_{i},\mathbf{\alpha}\rangle\not =0$ mod $p$. In particular,
if the support of $\mu$ is not parallel to a proper subspace of $\mathbf{Q}%
^{k}$, we have $\mathbf{y}_{1},...,\mathbf{y}_{k}\in V$, $\langle
\mathbf{y}_{i},\mathbf{\alpha}\rangle\not =0$ mod $p$, for some $i\in
\{1,...,k\}$.\smallskip

Henceforth, we will indicate by $B$ the matrix $(\mathbf{y}_{1}...\mathbf{y}%
_{k})$, where the vectors $\mathbf{y}_{1},...,\mathbf{y}_{k}$ are defined by
Lemma 2.4.\smallskip

\textbf{Lemma 2.5.}\ $\forall~e,j\in\mathbf{N}$, we have:
\[
{^{t}A}^{e}=\prod_{i=1}^{e}({^{t}A}-\lambda_{i}I)+\sum_{s=0}^{e-1}\left(
\lambda_{s+1}{^{t}A}^{e-s-1}\prod_{i=1}^{s}({^{t}A}-\lambda_{i}I)\right)
.\leqno1)
\]%
\[
{^{t}A}^{j}\prod_{i=1}^{e}({^{t}A}-\lambda_{i}I)=\sum_{h=e+1}^{d}%
\sum_{\substack{k_{1},...,k_{d-h+1}:\\k_{m}\geq0,\forall m=1,...,d-h+1,\\\sum
_{m=1}^{d-h+1}k_{m}=j-d+h}}\prod_{m=h}^{d}\lambda_{m}^{k_{d-m+1}}\prod
_{n=h+1}^{d}({^{t}A}-\lambda_{n}I)\prod_{i=1}^{e}({^{t}A}-\lambda_{i}%
I).\leqno2)
\]

\textbf{Proof.}

1) The proof is equal to the proof of Lemma 3.2 in \cite{Asci}.

2) We can suppose $e\leq d-1$, since otherwise the two members of 2) are equal
to the null matrix. Set:%
\begin{gather*}
H_{h,j}\equiv\left\{  k_{1},...,k_{d-h+1}:k_{m}\geq0,\forall
m=1,...,d-h+1,\sum_{m=1}^{d-h+1}k_{m}=j-d+h\right\}  ,\\
f(\lambda_{h},...,\lambda_{d},j-d+h)\equiv\sum_{\substack{H_{h,j}}}\prod
_{m=h}^{d}\lambda_{m}^{k_{d-m+1}}.
\end{gather*}

Proceed by induction on $j$; if $j=0$, $\forall~e=0,1,...,d-1$, the thesis is true.

Suppose that the thesis is true for $j$; then, for $j+1$,$\ \forall
\,e=0,1,...,d-1$:%
\begin{gather*}
{^{t}A}^{j+1}\prod_{i=1}^{e}({^{t}A}-\lambda_{i}I)={^{t}A}\left(  {^{t}A}%
^{j}\prod_{i=1}^{e}({^{t}A}-\lambda_{i}I)\right) \\
=\sum_{h=e+1}^{d}({^{t}A}-\lambda_{h}I+\lambda_{h}I)f(\lambda_{h}%
,...,\lambda_{d},j-d+h)\prod_{n=h+1}^{d}({^{t}A}-\lambda_{n}I)\prod_{i=1}%
^{e}({^{t}A}-\lambda_{i}I)\\
=({^{t}A}-\lambda_{e+1}I)f(\lambda_{e+1},...,\lambda_{d},j-d+e+1)\prod
_{n=e+2}^{d}({^{t}A}-\lambda_{n}I)\prod_{i=1}^{e}({^{t}A}-\lambda_{i}I)\\
+\sum_{h=e+1}^{d-1}\left(  \lambda_{h}f(\lambda_{h},...,\lambda_{d}%
,j-d+h)\prod_{n=h+1}^{d}({^{t}A}-\lambda_{n}I)\right. \\
\left.  +({^{t}A}-\lambda_{h+1}I)f(\lambda_{h+1},...,\lambda_{d}%
,j+1-d+h)\prod_{n=h+2}^{d}({^{t}A}-\lambda_{n}I)\right)  \prod_{i=1}^{e}%
({^{t}A}-\lambda_{i}I)\\
+\lambda_{d}\lambda_{d}^{j}\prod_{i=1}^{e}({^{t}A}-\lambda_{i}I).
\end{gather*}

Observe that%
\[
\lambda_{h}f(\lambda_{h},...,\lambda_{d},j-d+h)+f(\lambda_{h+1},...,\lambda
_{d},j+1-d+h)=f(\lambda_{h},...,\lambda_{d},j+1-d+h),
\]
since $f(\lambda_{h},...,\lambda_{d},j+1-d+h)$, a homogeneous polynomial of
degree $j+1-d+h$ in the variables $\lambda_{h},...,\lambda_{d}$, can be
obtained by multiplying $f(\lambda_{h},...,\lambda_{d},j-d+h)$, a homogeneous
polynomial of degree $j-d+h$ in the variables $\lambda_{h},...,\lambda_{d}$,
by the variable $\lambda_{h}$ and by summing up $f(\lambda_{h+1}%
,...,\lambda_{d},j+1-d+h)$, a homogeneous polynomial of degree $j+1-d+h$ in
the variables $\lambda_{h+1},...,\lambda_{d}$. Then:%
\begin{gather*}
{^{t}A}^{j+1}\prod_{i=1}^{e}({^{t}A}-\lambda_{i}I)=f(\lambda_{e+1}%
,...,\lambda_{d},j-d+e+1)\prod_{i=1}^{d}({^{t}A}-\lambda_{i}I)\\
+\sum_{h=e+1}^{d-1}f(\lambda_{h},...,\lambda_{d},j+1-d+h)\prod_{n=h+1}%
^{d}({^{t}A}-\lambda_{n}I)\prod_{i=1}^{e}({^{t}A}-\lambda_{i}I)+\lambda
_{d}^{j+1}\prod_{i=1}^{e}({^{t}A}-\lambda_{i}I)\\
=\sum_{h=e+1}^{d}\sum_{\substack{H_{h,j+1}}}\prod_{m=h}^{d}\lambda
_{m}^{k_{d-m+1}}\prod_{n=h+1}^{d}({^{t}A}-\lambda_{n}I)\prod_{i=1}^{e}({^{t}%
A}-\lambda_{i}I),
\end{gather*}
since $\prod_{i=1}^{d}({^{t}A}-\lambda_{i}I)$ is equal to the null matrix.
$\Box\medskip$

\section{Main results}

\textbf{Theorem 3.1.}\ Assume that $A$ has eigenvalues $\lambda_{1}%
,...,\lambda_{k}\in\mathbf{C}^{\ast}$, $|\lambda_{i}|\neq1$, for $i=1,...,k$,
and assume that the support of $\mu$ is not parallel to a proper subspace of
$\mathbf{Q}^{k}$ invariant under $A$. Then, there exists $c\in\mathbf{R}^{+}$
such that, for all $p\in\mathbf{N}$ such that $\gcd($det$(A),p)=$ $\gcd
($det$(B),p)=1$, and for all $n\geq c(\ln p)^{2}$, we have:
\[
\Vert P_{n}-U\Vert\leq\varepsilon(p)\text{,\qquad where }\lim_{p\rightarrow
\infty}\varepsilon(p)=0.
\]

\textbf{Proof.}\ From the lemmas 2.1 and 2.2:
\begin{equation}
\Vert P_{n}-U\Vert^{2}\leq{\frac{1}{4}}\sum_{\mathbf{\alpha}\in\mathbf{Z}%
_{p}^{k}-\{\mathbf{0}\}}|\widehat{{P}}_{n}(\mathbf{\alpha})|^{2};
\label{form1}%
\end{equation}%
\begin{equation}
|\widehat{{P}}_{n}(\mathbf{\alpha})|^{2}=\prod_{j=0}^{n-1}\left|
\widehat{{\mu}}\left(  {}^{t}A^{j}\mathbf{\alpha}\right)  \right|  ^{2}%
=\prod_{j=0}^{n-1}\left(  {\sum\limits_{\mathbf{u},\mathbf{v}\in\mathbf{Z}%
^{k}}{\mu(\mathbf{u})\mu(\mathbf{v})\cos\left(  {\frac{{2\pi}}{p}}\left\langle
\mathbf{u-v},{}^{t}A^{j}\mathbf{\alpha}\right\rangle \right)  }}\right)  .
\label{form2}%
\end{equation}

In order to estimate $\displaystyle\prod_{j=0}^{n-1}\left|  \widehat{{\mu}%
}\left(  {}^{t}A^{j}\mathbf{\alpha}\right)  \right|  ^{2}$, $\forall$
$\delta\in(0,\frac{1}{2})$ set:%
\[
L_{\delta}=\underset{i=1}{\overset{k}{\bigcup}}\left(  [0,1]^{i-1}%
\times\left[  \delta,1-\delta\right]  \times\lbrack0,1]^{k-i}\right)
,\quad\mathbf{\xi}_{j}={{\frac{{{^{t}A^{j}}}\mathbf{\alpha}}{p},}}%
\]
and indicate by $\left\{  \mathbf{\xi}_{j}\right\}  $ the vector whose
components are the fractional parts of the components of $\mathbf{\xi}_{j}$.
Consider the vectors $\mathbf{y}_{1},...,\mathbf{y}_{k}$ defined by Lemma 2.4;
then, $\forall\,m=1,...,k$,
\[
\mathbf{y}_{m}=A^{z_{m}}\mathbf{x}_{m}\text{,\quad where }\mathbf{x}_{m}\in
V,\ z_{m}\in\{0,1,...,d-1\}.
\]
Finally set $z=\underset{m=1,...,k}{\max}z_{m}$. The following results hold:\medskip

\textbf{Lemma 3.2.} There exists $b\in(0,1)$ such that, if $\left\{
\mathbf{\xi}_{i}\right\}  \in L_{\delta}$, then $\displaystyle\prod
_{j=0}^{i+z}\left|  \widehat{{\mu}}\left(  {}^{t}A^{j}\mathbf{\alpha}\right)
\right|  ^{2}\leq b$.

\textbf{Proof.} Consider the function $g:[0,1]^{k}\longrightarrow\lbrack0,1]$
defined by:
\[
g(\mathbf{t})=\prod_{j=0}^{z}\left|  \widehat{{\mu}}\left(  {}^{t}%
A^{j}p\mathbf{t}\right)  \right|  ^{2}=\prod_{j=0}^{z}\left(  {\sum
\limits_{\mathbf{u},\mathbf{v}\in\mathbf{Z}^{k}}{\mu(\mathbf{u})\mu
(\mathbf{v})\cos\left(  {2\pi}\left\langle \mathbf{u-v},{}^{t}A^{j}%
\mathbf{t}\right\rangle \right)  }}\right)  .
\]
If $t_{l}\in(0,1)$ for some $l\in\{1,...,k\}$, then $p\mathbf{t}\neq
\mathbf{0}$ (mod $p$) and, by Lemma 2.4,\ $\exists{}~m\in\{1,...,k\}$ such
that $\left\langle \mathbf{x}_{m},{}^{t}A^{z_{m}}p\mathbf{t}\right\rangle
\neq\mathbf{0}$ (mod $p$); then $\left\langle \mathbf{x}_{m},{}^{t}A^{z_{m}%
}\mathbf{t}\right\rangle \notin\mathbf{Z}$ and, by definition of $g$, we have
$g(\mathbf{t})<1$. Since $g$ is continuous and $L_{\delta}$ is closed and
bounded, then $g\left(  \mathbf{t}\right)  \leq b<1$, $\forall~\mathbf{t}\in
L_{\delta}$; in particular, if $\left\{  \mathbf{\xi}_{i}\right\}  \in
L_{\delta}$:%
\[
\prod_{j=0}^{i+z}\left|  \widehat{{\mu}}\left(  {}^{t}A^{j}\mathbf{\alpha
}\right)  \right|  ^{2}=\prod_{j=0}^{i-1}\left|  \widehat{{\mu}}\left(  {}%
^{t}A^{j}\mathbf{\alpha}\right)  \right|  ^{2}\prod_{l=0}^{z}\left|
\widehat{{\mu}}\left(  {}^{t}A^{l+i}\mathbf{\alpha}\right)  \right|  ^{2}\leq
g\left(  \left\{  \mathbf{\xi}_{i}\right\}  \right)  \leq b.~\Box
\]
\medskip

\textbf{Lemma 3.3. }There exist $\delta,\overline{c}\in\mathbf{R}^{+}$,
$\delta\in(0,\frac{1}{2})$, and $\overline{j}\in\mathbf{N}^{\ast}$,
$\overline{j}\leq\overline{c}\ln p$, such that, for any $p\in\mathbf{N}%
$\ sufficiently large and for any $\mathbf{\alpha}\in\mathbf{Z}_{p}%
^{k}-\{\mathbf{0}\}$, $\left\{  \mathbf{\xi}_{\overline{j}}\right\}  $ has a
component in $\left[  \delta,1-\delta\right]  $.

\textbf{Proof. }The proof follows from \ Lemma 3 in \cite{Hildebrand-07}.
$\Box\medskip$

Let $\displaystyle c>-\frac{k\overline{c}}{\ln b}$ and suppose $n\geq c(\ln
p)^{2}$; then, for $p$ sufficiently large, we have $n\geq rt$, where%
\[
t=\left\lfloor \overline{c}\ln p\right\rfloor +d,\quad\quad r=\left\lfloor
\overline{\overline{c}}\ln p\right\rfloor ,
\]
for some $\displaystyle\overline{\overline{c}}>-\frac{k}{\ln b}.$

Let $\overline{j}\in\mathbf{N}^{\ast}$ defined by Lemma 3.3. If
$\mathbf{\alpha}\in\widetilde{\mathbf{Z}}_{p}^{k}-\{\mathbf{0}\}$, from Lemma
3.2 and Lemma 3.3, since $\left\{  \mathbf{\xi}_{\overline{j}}\right\}  \in
L_{\delta}$:
\[
\prod_{j=0}^{\overline{j}+z}\left|  \widehat{{\mu}}\left(  {}^{t}%
A^{j}\mathbf{\alpha}\right)  \right|  ^{2}\leq b<1.
\]
Since $\overline{j}+z\leq t-1$, from (\ref{form2}):%
\[
|\widehat{{P}}_{t}(\mathbf{\alpha})|^{2}=\prod_{j=0}^{t-1}\left|
\widehat{{\mu}}\left(  {}^{t}A^{j}\mathbf{\alpha}\right)  \right|  ^{2}\leq
b.
\]
By repeating the previous arguments with ${}^{t}A^{it}\mathbf{\alpha}$ instead
of $\mathbf{\alpha}$, $\forall~i=0,1,...,r-1$, we have:%
\[
|\widehat{{P}}_{rt}(\mathbf{\alpha})|^{2}=\prod_{i=0}^{r-1}\prod_{j=0}%
^{t-1}\left|  \widehat{{\mu}}\left(  {}^{t}A^{j+it}\mathbf{\alpha}\right)
\right|  ^{2}\leq b^{r}.
\]
Then, from (\ref{form2}):%
\begin{gather*}
\sum_{\mathbf{\alpha}\in\mathbf{Z}_{p}^{k}-\{\mathbf{0}\}}|\widehat{{P}}%
_{n}(\mathbf{\alpha})|^{2}\leq\sum_{\mathbf{\alpha}\in\mathbf{Z}_{p}%
^{k}-\{\mathbf{0}\}}|\widehat{{P}}_{rt}(\mathbf{\alpha})|^{2}\leq b^{r}p^{k}\\
\leq\frac{b^{\overline{\overline{c}}\ln p}p^{k}}{b}=\frac{\exp\left(  \left(
\overline{\overline{c}}\ln b{+k}\right)  \ln p\right)  }{b}.
\end{gather*}
Since $\displaystyle\underset{p\longrightarrow\infty}{\lim}\frac{\exp\left(
\left(  \overline{\overline{c}}\ln b{+k}\right)  \ln p\right)  }{b}=0$, by
definition of $\overline{\overline{c}}$, from (\ref{form1}) we obtain the
statement of the theorem. $\Box\medskip$

\textbf{Theorem 3.4.}\ Assume that $A$ has eigenvalues $\lambda_{1}%
,...,\lambda_{k}\in\mathbf{C}^{\ast}$ such that $\lambda_{i}^{l_{i}}%
\in\mathbf{N}^{\ast}-\{1\}$, for some $l_{i}\in\mathbf{N}^{\ast}$, for any
$i=1,...,k$, and assume that the support of $\mu$ is not parallel to a proper
subspace of $\mathbf{Q}^{k}$ invariant under $A$. Then, there exists
$c\in\mathbf{R}^{+}$ such that, for all $p\in\mathbf{N}$ such that $\gcd
($det$(A),p)=$ $\gcd($det$(B),p)=1$, and for all $n\geq c\ln p\ln\ln p$, we
have:
\[
\Vert P_{n}-U\Vert\leq\varepsilon(p)\text{,\qquad where }\lim_{p\rightarrow
\infty}\varepsilon(p)=0.
\]

\textbf{Proof. }Set%
\begin{equation}
l=\operatorname{lcm}(l_{i}:i=1,...,k),~C=A^{l},~\sigma_{i}=\lambda_{i}%
^{l},~\forall~i=1,...,k. \label{def1}%
\end{equation}
Since $\overset{d}{\underset{i=1}{\prod}}(^{t}A-\lambda_{i}I)=0$, we have
$\overset{d}{\underset{i=1}{\prod}}(^{t}C-\sigma_{i}I)=0$. Moreover,
$\forall~E\subsetneqq\{1,...,d\}$, define:
\[
Y_{E}=\left\{  \mathbf{\alpha}\in\mathbf{Z}_{p}^{k}-\{\mathbf{0}%
\}:{\prod\limits_{i\in E}{\left(  ^{t}C{-\sigma_{i}I}\right)  }}%
\mathbf{\alpha}\neq\mathbf{0}\text{ (mod }p\text{)},\right.
\]%
\[
\left.  {\prod\limits_{i\in T}{\left(  ^{t}C{-\sigma_{i}I}\right)  }%
}\mathbf{\alpha}=\mathbf{0}\text{ (mod }p\text{)},\,\,\forall~T\subset
\{1,...,d\}:|E|+1\leq|T|\leq d\right\}  .
\]

The following relation holds:
\begin{equation}
\mathbf{Z}_{p}^{k}-\{\mathbf{0}\}=\bigcup\limits_{E\subsetneqq\{1,...,d\}}%
Y_{E}\mathbf{.} \label{rel2}%
\end{equation}
In fact, observe that, $\forall~\mathbf{\alpha}\in\mathbf{Z}_{p}%
^{k}-\{\mathbf{0}\}$, we can define:%
\begin{align*}
e  &  =\max\left\{  t\in\{0,1,...,d-1\}:{\prod\limits_{i\in E}{\left(
^{t}C{-\sigma_{i}I}\right)  }}\mathbf{\alpha}\neq\mathbf{0}\text{ (mod
}p\text{)},\right. \\
&  \left.  \text{for some }E\subsetneqq\{1,...,d\}\text{ such that
}|E|=t\right\}  ;
\end{align*}
then $\mathbf{\alpha}\in Y_{E}$, for some $E\subsetneqq\{1,...,d\}$ such that
$|E|=e$, and this implies (\ref{rel2}).

Moreover:%
\[
Y_{E}=\bigcup\limits_{\emptyset\neq S\subset\{1,...,k\}}Y_{S,E}\mathbf{,}%
\]
where, if the vectors $\mathbf{y}_{1},...,\mathbf{y}_{k}$ are defined by Lemma
2.4:%
\[
Y_{S,E}=\left\{  \mathbf{\alpha}\in Y_{E}:\left\langle \mathbf{y}_{m}%
,{\prod\limits_{i\in E}{\left(  ^{t}C{-\sigma_{i}I}\right)  }}\mathbf{\alpha
}\right\rangle \neq0\text{ (mod }p\text{)},~\forall~m\in S,\right.
\]%
\[
\left.  \left\langle \mathbf{y}_{m},{\prod\limits_{i\in E}{\left(
^{t}C{-\sigma_{i}I}\right)  }}\mathbf{\alpha}\right\rangle =0\text{ (mod
}p\text{)},~\forall~m\notin S\right\}  .
\]
Then:%
\begin{equation}
\sum_{\mathbf{\alpha}\in\mathbf{Z}_{p}^{k}-\{\mathbf{0}\}}|\widehat{{P}}%
_{n}(\mathbf{\alpha})|^{2}=\sum_{E\subsetneqq\{1,...,d\}}\sum_{\emptyset\neq
S\subset\{1,...,k\}}\sum_{Y_{S,E}}|\widehat{{P}}_{n}(\mathbf{\alpha})|^{2}.
\label{form3}%
\end{equation}

If $\emptyset\neq S\subset\{1,...,k\}$, reorder the set $S$ in the following
manner:%
\[
S=\left\{  m_{1,S},...,m_{|S|,S}\right\}  ,\quad\text{where }m_{i,S}%
<m_{j,S}\Leftrightarrow i<j.
\]
Then, $\forall~h=1,...,|S|$:
\[
\mathbf{y}_{m_{h,S}}=A^{z_{m_{h,S}}}\mathbf{x}_{m_{h,S}}\text{,\quad where
}\mathbf{x}_{m_{h,S}}\in V,\ z_{m_{h,S}}\in\{0,1,...,d-1\}.
\]
Set $\overline{\mathbf{x}}_{h,S}\equiv\mathbf{x}_{m_{h,S}}$, $\overline
{\mathbf{y}}_{h,S}\equiv\mathbf{y}_{m_{h,S}}$, $\overline{z}_{h,S}\equiv
z_{m_{h,S}}$ and let $\overline{\mathbf{u}}_{h,S}$, $\overline{\mathbf{v}%
}_{h,S}$ the vectors of the support of $\mu$ such that $\overline{\mathbf{u}%
}_{h,S}-\overline{\mathbf{v}}_{h,S}=\overline{\mathbf{x}}_{h,S}$. From Lemma
2.2, $\forall~\mathbf{\alpha}\in\mathbf{Z}_{p}^{k}-\{\mathbf{0}\}$, we have:%
\begin{gather}
|\widehat{{P}}_{n}(\mathbf{\alpha})|^{2}\leq\prod_{j=0}^{\left\lfloor
\frac{n-1}{l}\right\rfloor }\left|  \widehat{{\mu}}\left(  {}^{t}%
C^{j}\mathbf{\alpha}\right)  \right|  ^{2}\nonumber\\
\leq\prod_{h=1}^{|S|}\prod_{j\in M_{h,S}}\left(  1-2{{\mu(\overline
{\mathbf{u}}_{h,S})\mu(\overline{\mathbf{v}}_{h,S})+2{\mu(\overline
{\mathbf{u}}_{h,S})\mu(\overline{\mathbf{v}}_{h,S})}\cos\left(  {\frac{{2\pi}%
}{p}}\left\langle \overline{\mathbf{y}}_{h,S},{}^{t}C^{j}\mathbf{\alpha
}\right\rangle \right)  }}\right)  , \label{formin}%
\end{gather}
where $\displaystyle M_{h,S}\equiv\left\{  (h-1)\left\lfloor \frac
{\left\lfloor \frac{n-1}{l}\right\rfloor }{|S|}\right\rfloor
,...,h\left\lfloor \frac{\left\lfloor \frac{n-1}{l}\right\rfloor }%
{|S|}\right\rfloor -\overline{z}_{h,S}-1\right\}  $, $\forall~h=1,...,|S|$.

If $\mathbf{\alpha}\in Y_{S,E}$, for some $\emptyset\neq S\subset\{1,...,k\}$
and $E\subsetneqq\{1,...,d\}$, reorder the numbers $\sigma_{1},...,\sigma_{d}$
so that the first $|E|$ correspond to the set $\{\sigma_{i}:i\in E\}$.
Moreover, $\forall~j,n\in\mathbf{N}$, $\forall~h=1,...,|S|$, set:%
\[
a_{h,S,E}=\left\langle \overline{\mathbf{y}}_{h,S},\overset{|E|}%
{{\prod\limits_{i=1}}}{{\left(  ^{t}C{-\sigma_{i}I}\right)  }}\mathbf{\alpha
}\right\rangle \in\mathbf{Z}_{p}-\{0\},\quad\xi_{j,n,h,S}=\left\{
\left\langle \overline{\mathbf{y}}_{h,S},^{t}C^{j}\overset{n}{{\prod
\limits_{i=1}}}{{\left(  ^{t}C{-\sigma_{i}I}\right)  }}\frac{\mathbf{\alpha}%
}{p}\right\rangle \right\}
\]
(the fractional part of $\displaystyle\left\langle \overline{\mathbf{y}}%
_{h,S},^{t}C^{j}\overset{n}{{\prod\limits_{i=1}}}{{\left(  ^{t}C{-\sigma_{i}%
I}\right)  }}\frac{\mathbf{\alpha}}{p}\right\rangle $). Observe that, from
Lemma 2.5, 1):%
\[
\xi_{j+|E|,0,h,S}=\left\{  \left\langle \overline{\mathbf{y}}_{h,S},^{t}%
C^{j}\overset{|E|}{{\prod\limits_{i=1}}}\left(  ^{t}C-\sigma_{i}{I}\right)
\frac{\mathbf{\alpha}}{p}\right\rangle +\sum_{m=0}^{|E|-1}\sigma_{m+1}%
\xi_{j+|E|-m-1,m,h,S}\right\}  .
\]
Moreover, use Lemma 2.5, 2) and the definition of $Y_{E}$; $\forall$
$\mathbf{\alpha}\in Y_{S,E}$, $\forall~j\in\mathbf{N}$, in the right member of
2) multiplied by $\mathbf{\alpha}$ (with $e=|E|$, $A=C$ and $\lambda
_{i}=\sigma_{i}$), only $\sigma_{d}^{j}\overset{|E|}{{\prod\limits_{i=1}}%
}{{\left(  ^{t}C{-\sigma_{i}I}\right)  }}\mathbf{\alpha}$ is different from
$0$ (corresponding to $h=m=d$, $k_{1}=j$). Then, $\forall~h=1,...,|S|$:%
\begin{gather}
^{t}C^{j}\overset{|E|}{{\prod\limits_{i=1}}}{{\left(  ^{t}C{-\sigma_{i}%
I}\right)  }}\mathbf{\alpha=}\sigma_{d}^{j}\overset{|E|}{{\prod\limits_{i=1}}%
}{{\left(  ^{t}C{-\sigma_{i}I}\right)  }}\mathbf{\alpha\Rightarrow
}\left\langle \overline{\mathbf{y}}_{h,S},^{t}C^{j}\overset{|E|}%
{{\prod\limits_{i=1}}}{{\left(  ^{t}C{-\sigma_{i}I}\right)  }}\frac
{\mathbf{\alpha}}{p}\right\rangle =\left.  \sigma_{d}^{j}\right.
\frac{a_{h,S,E}}{p}\nonumber\\
\Rightarrow\xi_{j+|E|,0,h,S}=\left\{  \left.  \sigma_{d}^{j}\right.
\frac{a_{h,S,E}}{p}+\sum_{m=0}^{|E|-1}\sigma_{m+1}\xi_{j+|E|-m-1,m,h,S}%
\right\}  .\label{csi-sig}%
\end{gather}

Let $g_{h,S}:[0,1]\longrightarrow\lbrack0,1]$ the function defined by:
\[
g_{h,S}(t)=1-2{{\mu(\overline{\mathbf{u}}_{h,S})\mu(\overline{\mathbf{v}%
}_{h,S})+2{\mu(\overline{\mathbf{u}}_{h,S})\mu(\overline{\mathbf{v}}_{h,S})}}%
}\cos(2\pi t).
\]
From (\ref{formin}), we have:%
\begin{gather}
\sum_{Y_{S,E}}|\widehat{{P}}_{n}(\mathbf{\alpha})|^{2}\leq\sum
_{\substack{a_{h,S,E}\in\mathbf{Z}_{p}-\{0\},\\\forall~h=1,...,|S|}%
}\prod_{h=1}^{|S|}\prod_{j\in M_{h,S}}g_{h,S}\left(  \xi_{j,0,h,S}\right)
\nonumber\\
=\prod_{h=1}^{|S|}\left.  \sum_{a_{h,S,E}\in\mathbf{Z}_{p}-\{0\}}\right.
\prod_{j\in M_{h,S}}g_{h,S}\left(  \xi_{j,0,h,S}\right)  . \label{sigma}%
\end{gather}

Set $\displaystyle L={\left[  {\frac{1}{{2^{d}\sigma^{d+1}}}},\ 1-{\frac
{1}{{2^{d}\sigma^{d+1}}}}\right]  }$, where $\sigma\equiv\underset
{i=1,...,k}{\max}\sigma_{i}$.

Observe that $g_{h,S}$ is continuous, $g_{h,S}(t)=1\Leftrightarrow
t\in\{0,1\}$; then, $g_{h,S}$ has a maximum $b_{h,S}<1$ in $L$, since $L$ is
closed and bounded; in particular:%
\begin{equation}
g_{h,S}\left(  \xi_{j,0,h,S}\right)  \leq\left\{
\begin{array}
[c]{ll}%
b_{h,S} & \text{if }\xi_{j,0,h,S}\in L\\
1 & \text{otherwise}.
\end{array}
\right.  \label{gi-bi}%
\end{equation}

The following result follows:\medskip

\textbf{Lemma 3.5. }Suppose $\xi_{j,0,h,S},\xi_{j+1,0,h,S},...,\xi
_{j+e-1,0,h,S}\not \in L$, for some $j\in\mathbf{N}^{\ast}$, $e\in\mathbf{N}$;
then, for any $s\in\{0,1,...,e-1\}$ and for any$\ r\in\{j,j+1,...,j+e-s-1\}$:
\begin{equation}
\xi_{r,s,h,S}\in\left[  -{\frac{1}{{2^{d-s}\sigma^{d-s+1}}},\frac{1}%
{{2^{d-s}\sigma^{d-s+1}}}}\right]  \text{ (mod }\mathbf{Z}\text{)}.
\label{beta3}%
\end{equation}
In particular:%

\begin{equation}
\xi_{j+e-s-1,s,h,S}\in\left[  -{\frac{1}{{2^{d-s}\sigma^{3}}},\frac
{1}{{2^{d-s}\sigma^{3}}}}\right]  \text{ (mod }\mathbf{Z}\text{)}.
\label{beta4}%
\end{equation}

\textbf{Proof.} Prove the lemma by induction on $s$; by hypothesis, if $s=0$
and $r\in\{j,j+1,...,j+e-1\}$, the thesis is true.

Suppose that the thesis is true for $s=n$; then, for $s=n+1,\,\forall
~r=j,j+1,...,j+e-n-2$, we have:
\[
{{^{t}C^{r}}}\prod\limits_{i=1}^{n+1}\left(  {{^{t}C}-\sigma_{i}I}\right)
\,={{^{t}C^{r+1}}}\prod\limits_{i=1}^{n}{}\left(  {{^{t}C}-\sigma_{i}%
I}\right)  -\sigma_{n+1}{{^{t}C^{r}}}\prod\limits_{i=1}^{n}{}\left(  {{^{t}%
C}-\sigma_{i}I}\right)  .
\]
By the inductive hypothesis:
\begin{align*}
\xi_{r,n+1,h,S}  &  \in2\sigma\left[  -{\frac{1}{{2^{d-n}\sigma^{d-n+1}}%
},\frac{1}{{2^{d-n}\sigma^{d-n+1}}}}\right]  \text{ (mod }\mathbf{Z}\text{)}\\
&  =\left[  -{\frac{1}{{2^{d-(n+1)}\sigma^{d-(n+1)+1}}},\frac{1}%
{{2^{d-(n+1)}\sigma^{d-(n+1)+1}}}}\right]  \text{ (mod }\mathbf{Z}\text{)}.
\end{align*}
Thus, we have (\ref{beta3}). In particular, since $d-s\geq2$, (\ref{beta4})
follows. $\Box$\medskip

In order to finish the proof of Theorem 3.4, we will borrow some arguments
from the papers \cite{Asci}, \cite{Chung} and \cite{Hildebrand-93}.

Fix $h,S,E$, set $a=a_{h,S,E}$ and consider the expansion of ${\frac{a}{p}}$
in base $\sigma_{d}$:
\[
{\frac{a}{p}}=\ 0.a_{1}a_{2}a_{3}...\quad\text{Define:}%
\]%
\[
t=\left\lceil \log_{\sigma_{d}}p\right\rceil ~~(\text{then, }\sigma_{d}%
^{t-1}<p<\sigma_{d}^{t}),\quad r=r_{h,S}=\left\lfloor \frac{|M_{h,S}|}%
{t}\right\rfloor =\left\lfloor \left.  \left(  \left\lfloor \frac{\left\lfloor
\frac{n-1}{l}\right\rfloor }{|S|}\right\rfloor -\overline{z}_{h,S}\right)
\right/  t\right\rfloor .
\]

Moreover, recall that a ''generalized alternation'' between two consecutive
digits $a_{j}a_{j+1}$ of the expansion is defined as either the case
$a_{j}\neq a_{j+1}$ or the case $a_{j}=a_{j+1}\notin\{0,\sigma_{d}-1\}$.\medskip

\textbf{Lemma 3.6. }Suppose that there is a generalized alternation between
the digits $j+1,j+2$ - th of the expansion of ${\frac{a}{p}}$; then
$\xi_{j+i,0,h,S}\in L$, for some $i\in\{0,1,...,|E|\}$.

\textbf{Proof. }The assumption imply:
\[
\left\{  \sigma{_{d}^{j}\frac{a}{p}}\right\}  \in\left[  {\frac{1}{\sigma
_{d}^{2}}},\ 1-{\frac{1}{\sigma_{d}^{2}}}\right]  \subset\left[  {\frac
{1}{\sigma^{2}}},\ 1-{\frac{1}{\sigma^{2}}}\right]  .
\]
If $\xi_{j,0,h,S},\xi_{j+1,0,h,S},...,\xi_{j+|E|-1,0,h,S}\not \in L$, then,
from (\ref{csi-sig}) and Lemma 3.5:
\begin{align*}
\xi_{j+|E|,0,h,S}  &  \in\left(  \left[  {\frac{1}{\sigma^{2}}},\ 1-{\frac
{1}{\sigma^{2}}}\right]  +\left[  -\sum_{m=0}^{|E|-1}{\frac{\sigma_{m+1}%
}{2^{d-m}\sigma^{3}}},\ \sum_{m=0}^{!E|-1}{\frac{\sigma_{m+1}}{2^{d-m}%
\sigma^{3}}}\right]  \right) \\
&  \subset\left(  \left[  {\frac{1}{\sigma^{2}}},\ 1-{\frac{1}{\sigma^{2}}%
}\right]  +\left[  -{\frac{1}{\sigma^{2}}}\sum_{m=0}^{d-1}{\frac{1}{2^{d-m}}%
},\ {\frac{1}{\sigma^{2}}}\sum_{m=0}^{d-1}{\frac{1}{2^{d-m}}}\right]  \right)
\\
&  =\left[  {\frac{1}{2^{d}\sigma^{2}}},\ 1-{\frac{1}{2^{d}\sigma^{2}}%
}\right]  \subset L.~\Box
\end{align*}
\medskip

Consider the first $rt$ integer numbers of the set $M_{h,S}$ and partition
such numbers into $r$ disjoint sets $M_{i}=M_{i,h,S},\ 1\leq i\leq r$, each of
length $t$, such that, if $i<j$, $x\in M_{i}$ and $y\in M_{j}$, then $x<y$.
Moreover, $\forall\ i=1,...,r$, consider the block of digits
\[
B_{a,i}=B_{a,i,h,S}=\{a_{j}:j\in M_{i}\}
\]
and, $\forall\ B\subset\{a_{1},a_{2},...\}$, $B$ made up of consecutive
digits, indicate by $A(B)$ the number of generalized alternations in $B$.
Finally, $\forall~D\subset\mathbf{N}$, set:
\[
C(D)=C(D,h,S)=|\{j\in D:\xi_{j,0,h,S}\in L\}|.
\]

Suppose $A(\{a_{j}:j\in M_{h,S}\})=md$, for some $m\in\mathbf{N}$. Since
$|E|\leq d-1$, from Lemma 3.6 we deduce $C(M_{h,S})\geq m$; in general:
\begin{equation}
C(M_{h,S})\geq\left\lfloor {\frac{A(\{a_{j}:j\in M_{h,S}\})}{d}}\right\rfloor
\geq\frac{{\sum\limits_{i=1}^{r}{A(B_{a,i})}}}{d}-1. \label{ci}%
\end{equation}

It is possible to prove the following two results:\medskip

\textbf{Lemma 3.7.}\ $\forall~i\in\{1,...,r\}$, as $a$ ranges in
$\mathbf{Z}_{p}-\{0\}$, the blocks $B_{a,i}$ are distinct and have at least
one generalized alternation. Moreover, $\forall~\,i,j\in\{1,...,r\}$:
\[
\{B_{a,i}:a\in\mathbf{Z}_{p}-\{0\}\}=\{B_{a,j}:a\in\mathbf{Z}_{p}-\{0\}\}.
\]

\textbf{Lemma 3.8.}%
\[
{\sum_{j=1}^{s}\prod_{i=1}^{r}a_{\pi_{i}(j)}\leq\sum_{j=1}^{s}a_{j}^{r}},
\]
where, $\forall~i=1,...,r$, $\forall~j=1,...,s$, $\pi_{i}$ is a permutation of
$\{1,...,s\}$ and$\ a_{j}\geq0$.\medskip

By utilizing the formulas (\ref{sigma}), (\ref{gi-bi}), (\ref{ci}) and the
lemmas 3.7 and 3.8, we have:
\begin{align}
\sum_{Y_{S,E}}|\widehat{{P}}_{n}(\mathbf{\alpha})|^{2}  &  \leq\prod
_{h=1}^{|S|}\left.  \sum_{a\in\mathbf{Z}_{p}-\{0\}}\right.  \prod_{j\in
M_{h,S}}g_{h,S}\left(  \xi_{j,0,h,S}\right)  \leq\prod_{h=1}^{|S|}\left.
\sum_{a\in\mathbf{Z}_{p}-\{0\}}\right.  b_{h,S}^{C(M_{h,S})}\nonumber\\
&  \leq\frac{1}{{b_{S}}}\prod_{h=1}^{|S|}\sum\limits_{a\in\mathbf{Z}%
_{p}-\{0\}}\prod_{i=1}^{r}f_{h,S}^{A(B_{a,i})}\quad(\text{where }b_{S}%
=\prod_{h=1}^{|S|}b_{h,S},~f_{h,S}=\sqrt[d]{b_{h,S}}<1)\nonumber\\
&  \leq\frac{1}{{b_{S}}}\prod_{h=1}^{|S|}\sum\limits_{a\in\mathbf{Z}%
_{p}-\{0\}}f_{h,S}^{rA(B_{a,1})}\leq\frac{1}{{b_{S}}}\prod_{h=1}^{|S|}%
\sum_{\substack{\text{length }B=t\\A(B)>0}}f_{h,S}^{rA(B)}. \label{mulmin}%
\end{align}

Indicate with $M(j)$ the number of blocks of length $t$ with $A(B)=j$; then:
\[
M(j)\leq\left(
\begin{array}
[c]{c}%
t-1\\
j
\end{array}
\right)  \sigma_{d}^{j+1}\leq\left(
\begin{array}
[c]{c}%
t\\
j
\end{array}
\right)  \sigma_{d}^{j+1}%
\]%
\begin{align}
&  \Rightarrow\sum_{\substack{\text{length }B=t\\A(B)>0}}f_{h,S}^{rA(B)}%
\leq\sum_{j=1}^{t-1}M(j)f_{h,S}^{rj}\leq\sigma_{d}\sum_{j=1}^{t}\left(
\begin{array}
[c]{c}%
t\\
j
\end{array}
\right)  \left(  \sigma_{d}{f_{h,S}^{r}}\right)  ^{j}\nonumber\\
&  =\sigma_{d}\left(  \left(  1+\sigma_{d}{f_{h,S}^{r}}\right)  ^{t}-1\right)
\leq\sigma_{d}\left(  \exp\left(  \sigma_{d}t{f_{h,S}^{r}}\right)  -1\right)
. \label{exp}%
\end{align}

Suppose $\displaystyle c_{h,S}>-\frac{l|S|}{\ln f_{h,S}\ln\sigma_{d}}$, $n\geq
c_{h,S}\ln p\ln\ln p$. Since $\displaystyle t=\left\lceil \frac{\ln p}%
{\ln\sigma_{d}}\right\rceil \leq\frac{\ln p}{\ln\sigma_{d}}+1$, $\displaystyle
\exists~\overline{c}_{h,S}\in\left(  \frac{l|S|}{\ln f_{h,S}\ln\sigma_{d}%
},c_{h,S}\right)  $ such that, for $p$ sufficiently large:
\begin{align*}
r  &  =\left\lfloor \left.  \left(  \left\lfloor \frac{\left\lfloor \frac
{n-1}{l}\right\rfloor }{|S|}\right\rfloor -\overline{z}_{h,S}\right)  \right/
t\right\rfloor \geq\frac{\ln\sigma_{d}\overline{c}_{h,S}}{l|S|}\ln\ln p,\quad
t\leq\frac{2\ln p}{\ln\sigma_{d}}\\
&  \Rightarrow\sigma_{d}t{f_{h,S}^{r}\leq}\frac{2\sigma_{d}}{\ln\sigma_{d}%
}\exp\left(  \left(  1+\frac{\ln f_{h,S}\ln\sigma_{d}\overline{c}_{h,S}}%
{l|S|}\right)  \ln\ln p\right)  \equiv\gamma_{h,S}(p),
\end{align*}
where $\displaystyle\underset{p\rightarrow\infty}{\lim}\gamma_{h,S}(p)=0$, by
definition of $\overline{c}_{h,S}$. Finally:
\[
(\ref{exp})\leq\sigma_{d}\left(  \exp\left(  \gamma_{h,S}(p)\right)
-1\right)  \equiv\varepsilon_{h,S}(p),\quad\text{where }\underset
{p\rightarrow\infty}{\lim}\varepsilon_{h,S}(p)=0.
\]
From the formulas (\ref{form1}), (\ref{form3}), (\ref{mulmin}) and
(\ref{exp}), we obtain the thesis of the theorem, with $c=\displaystyle
\underset{h,S}{\max}{c_{h,S}}$.$\ \Box$\medskip

\textbf{Theorem 3.9.}\ Assume that $A$ has eigenvalues $\lambda_{1}%
,...,\lambda_{d},...,\lambda_{k}\in\mathbf{C}^{\ast}$ such that $\lambda
_{i}^{l_{i}}\in\mathbf{N}^{\ast}$, for some $l_{i}\in\mathbf{N}^{\ast}$, for
any $i=1,...,k$, where $|\lambda_{1}|=1$\ and$\ |\lambda_{i}|>1\ \forall
~i=2,...,d$, and assume that the support of $\mu$ is not parallel to a proper
subspace of $\mathbf{Q}^{k}$ invariant under $A$. Then, there exist
$\alpha,c\in\mathbf{R}^{+}$ and $N\in\mathbf{N}$ such that, for all
$p\in\mathbf{N}$ such that $p>N,\ \gcd(\det(A),p)=\gcd(\det(B),p)=1$, and for
all $n\geq cp^{2}$, we have:
\[
\Vert P_{n}-U\Vert\leq{2}^{k-1}\exp\left(  -{\frac{\alpha n}{p^{2}}}\right)
.
\]

\textbf{Proof.} Define $l$, $C$ and $\sigma_{i}$, $\forall~i=1,...,k$, as in
(\ref{def1}); then, $\overset{d}{\underset{i=1}{\prod}}(^{t}C-\sigma_{i}I)=0$.
Suppose $d>1$ and define, $\forall~\{1\}\subset E\subsetneqq\{1,...,d\}$:%
\[
Z=\left\{  \mathbf{\alpha}\in\mathbf{Z}_{p}^{k}-\{\mathbf{0}\}:{{\left(
^{t}C{-I}\right)  }}\mathbf{\alpha}\neq\mathbf{0}\text{ (mod }p\text{)}%
\right\}  ,
\]%
\begin{align*}
Z_{E}  &  =\left\{  \mathbf{\alpha}\in Z:{\prod\limits_{i\in E}{\left(
^{t}C{-\sigma_{i}I}\right)  }}\mathbf{\alpha}\neq\mathbf{0}\text{ (mod
}p\text{)},\right. \\
&  \left.  {\prod\limits_{i\in T}{\left(  ^{t}C{-\sigma_{i}I}\right)  }%
}\mathbf{\alpha}=\mathbf{0}\text{ (mod }p\text{)},\,\,\forall~T\subset
\{1,...,d\}:|E|+1\leq|T|\leq d\right\}  .
\end{align*}
The following relation is analogous to (\ref{rel2}):
\begin{equation}
Z=\bigcup\limits_{\{1\}\subset E\subsetneqq\{1,...,d\}}Z_{E}\mathbf{.}
\label{rel3}%
\end{equation}
Then:%
\begin{equation}
\sum_{\mathbf{\alpha}\in\mathbf{Z}_{p}^{k}-\{\mathbf{0}\}}|\widehat{{P}}%
_{n}(\mathbf{\alpha})|^{2}=\sum\limits_{\{1\}\subset E\subsetneqq
\{1,...,d\}}\sum_{Z_{E}}|\widehat{{P}}_{n}(\mathbf{\alpha})|^{2}+\sum_{Z^{c}%
}|\widehat{{P}}_{n}(\mathbf{\alpha})|^{2}. \label{form4}%
\end{equation}
Analogously to the proof of Theorem 3.4, where we have valued $\displaystyle
\sum_{Y_{S,E}}|\widehat{{P}}_{n}(\mathbf{\alpha})|^{2}$, by utilizing the fact
that $\sigma_{i}>1$, $\forall~i=2,...,d$, we can prove that $\exists
~c_{1},c_{2}\in\mathbf{R}^{+}$ such that:
\begin{equation}
\sum\limits_{\{1\}\subset E\subsetneqq\{1,...,d\}}\sum_{Z_{E}}|\widehat{{P}%
}_{n}(\mathbf{\alpha})|^{2}\leq\varepsilon_{1}(p,n),\quad\text{where
}\varepsilon_{1}(p,n)=c_{1}\ln p\exp\left(  -\frac{n}{c_{2}\ln p}\right)  .
\label{eps1}%
\end{equation}
Indeed, we must estimate only the sum $\displaystyle\sum\limits_{Z^{c}%
}|\widehat{{P}}_{n}(\mathbf{\alpha})|^{2}$. Observe that%
\[
Z^{c}=\bigcup\limits_{\emptyset\neq S\subset\{1,...,k\}}\overline{Z}%
_{S}\mathbf{,}%
\]
where, if the vectors $\mathbf{y}_{1},...,\mathbf{y}_{k}$ are defined by Lemma
2.4:%
\[
\overline{Z}_{S}=\left\{  \mathbf{\alpha}\in Z^{c}:\left\langle \mathbf{y}%
_{m},\mathbf{\alpha}\right\rangle \neq0\text{ (mod }p\text{)},~\forall~m\in
S,\right.
\]%
\[
\left.  \left\langle \mathbf{y}_{m},\mathbf{\alpha}\right\rangle =0\text{ (mod
}p\text{)},~\forall~m\notin S\right\}  .
\]
Then:%
\begin{equation}
\sum_{Z^{c}}|\widehat{{P}}_{n}(\mathbf{\alpha})|^{2}=\sum_{\emptyset\neq
S\subset\{1,...,k\}}\sum_{\overline{Z}_{S}}|\widehat{{P}}_{n}(\mathbf{\alpha
})|^{2}. \label{form5}%
\end{equation}

If$~\emptyset\neq S\subset\{1,...,k\}$ and $h\in\{1,...,|S|\}$,\ define
$\overline{\mathbf{x}}_{h,S}$, $\overline{\mathbf{y}}_{h,S}$, $\overline
{z}_{h,S}$, $\overline{\mathbf{u}}_{h,S}$ and $\overline{\mathbf{v}}_{h,S}$ as
in the proof of Theorem 3.4. Moreover, set:%
\[
a_{h,S}=\left\langle \overline{\mathbf{y}}_{h,S},\mathbf{\alpha}\right\rangle
.
\]
Observe that, $\forall~j\in\mathbf{N}$ and$\,\forall~\mathbf{\alpha}\in Z^{c}%
$, we have:
\[
\left\langle \overline{\mathbf{y}}_{h,S},{{}^{t}C^{j}}\mathbf{\alpha
}\right\rangle =a_{h,S}~\text{(mod }p\text{)}.
\]
Then, from the formula (\ref{formin}), $\forall~\mathbf{\alpha}\in Z^{c}$:%
\begin{gather}
|\widehat{{P}}_{n}(\mathbf{\alpha})|^{2}\leq\prod_{j=0}^{\left\lfloor
\frac{n-1}{l}\right\rfloor }\left|  \widehat{{\mu}}\left(  {}^{t}%
C^{j}\mathbf{\alpha}\right)  \right|  ^{2}\nonumber\\
\leq\prod_{h=1}^{|S|}\left(  1-2{{\mu(\overline{\mathbf{u}}_{h,S}%
)\mu(\overline{\mathbf{v}}_{h,S})+2{\mu(\overline{\mathbf{u}}_{h,S}%
)\mu(\overline{\mathbf{v}}_{h,S})}\cos\left(  {\frac{{2\pi}}{p}}%
a_{h,S}\right)  }}\right)  ^{|M_{h,S}|},\label{formin2}%
\end{gather}
where $\displaystyle|M_{h,S}|=\left\lfloor \frac{\left\lfloor \frac{n-1}%
{l}\right\rfloor }{|S|}\right\rfloor -\overline{z}_{h,S}\geq\frac{n-1}%
{l|S|}-1-(d-1)\geq\frac{n-(dlk+1)}{lk}$. Then:
\begin{align}
&  \sum_{\overline{Z}_{S}}|\widehat{{P}}_{n}(\mathbf{\alpha})|^{2}\nonumber\\
&  \leq\sum_{\substack{a_{h,S}\in\mathbf{Z}_{p}-\{0\},\\\forall~h=1,...,|S|}%
}\prod_{h=1}^{|S|}\left(  1-2{{\mu(\overline{\mathbf{u}}_{h,S})\mu
(\overline{\mathbf{v}}_{h,S})+2{\mu(\overline{\mathbf{u}}_{h,S})\mu
(\overline{\mathbf{v}}_{h,S})}\cos\left(  {\frac{{2\pi}}{p}}a_{h,S}\right)  }%
}\right)  ^{|M_{h,S}|}\nonumber\\
&  =\prod_{h=1}^{|S|}\sum_{a_{h,S}\in\mathbf{Z}_{p}-\{0\}}\left(
1-2{{\mu(\overline{\mathbf{u}}_{h,S})\mu(\overline{\mathbf{v}}_{h,S}%
)+2{\mu(\overline{\mathbf{u}}_{h,S})\mu(\overline{\mathbf{v}}_{h,S})}%
\cos\left(  {\frac{{2\pi}}{p}}a_{h,S}\right)  }}\right)  ^{|M_{h,S}%
|}.\label{zeta}%
\end{align}
Note that $\displaystyle-1+\cos x\leq-{\frac{2}{\pi^{2}}}x^{2},\ \forall
~x\in\lbrack-\pi,\pi]$. Furthermore, if $a_{h,S}\in\mathbf{Z}_{p}-\{0\}$, we
can suppose:
\[
a_{h,S}\in\mathbf{Z}^{\ast}\cap\left[  {-\frac{p-1}{2},\frac{p}{2}}\right]
\Rightarrow{\frac{2\pi}{p}}a_{h,S}\in\lbrack-\pi,\pi].
\]
Then, $\forall~h=1,...,|S|$:
\begin{align}
&  \sum_{a_{h,S}\in\mathbf{Z}_{p}-\{0\}}\left(  1-2{{\mu(\overline{\mathbf{u}%
}_{h,S})\mu(\overline{\mathbf{v}}_{h,S})+2{\mu(\overline{\mathbf{u}}_{h,S}%
)\mu(\overline{\mathbf{v}}_{h,S})}\cos\left(  {\frac{{2\pi}}{p}}%
a_{h,S}\right)  }}\right)  ^{|M_{h,S}|}\nonumber\\
&  \leq2{\sum\limits_{a_{h,S}\in\left(  \mathbf{Z}\cap\left[  1{,\frac{p}{2}%
}\right]  \right)  }\left(  {1-\frac{{16}}{{p^{2}}}{\mu(\overline{\mathbf{u}%
}_{h,S})\mu(\overline{\mathbf{v}}_{h,S})}a_{h,S}^{2}}\right)  ^{{|M}_{h,S}|}%
}\nonumber\\
&  \leq2{\sum\limits_{a_{h,S}\in\mathbf{N}^{\ast}}\exp\left(  {-\frac{{16}%
}{{lkp^{2}}}{\mu(\overline{\mathbf{u}}_{h,S})\mu(\overline{\mathbf{v}}_{h,S}%
)}(n-(dlk+1))a_{h,S}^{2}}\right)  }.\label{sigma2}%
\end{align}

Let $\overline{\tau}\in(0,1)$ such that $2\overline{\tau}^{3}+\overline{\tau
}^{2}-1=0$ ($\displaystyle\Leftrightarrow\frac{\overline{\tau}^{3}%
}{1-\overline{\tau}^{2}}=\frac{1}{2}$) and set:%
\[
t_{h,S}={\exp\left(  {-\frac{{16}}{{lkp^{2}}}{\mu(\overline{\mathbf{u}}%
_{h,S})\mu(\overline{\mathbf{v}}_{h,S})}(n-(dlk+1))}\right)  ,\quad}%
\overline{c}=\underset{h,S}{\max}{\frac{-lk{\ln}\overline{\tau}}%
{{16{\mu(\overline{\mathbf{u}}_{h,S})\mu(\overline{\mathbf{v}}_{h,S})}}}.}%
\]
Let $c>\overline{c}$, $p$ sufficiently large and $n\geq cp^{2}$; then:%
\[
{n-(dlk+1)}{\geq}\overline{c}p^{2}\Rightarrow t_{h,S}\leq\overline{\tau},
\]
from which
\begin{align*}
(\ref{sigma2})  &  =2\left(  {\sum\limits_{a_{h,S}\in\mathbf{N}^{\ast}}%
}t_{h,S}^{a_{h,S}^{2}}\right)  =2\left(  t_{h,S}{+\sum\limits_{a_{h,S}\geq2}%
}t_{h,S}^{a_{h,S}^{2}}\right)  \leq2\left(  t_{h,S}{+\sum\limits_{a_{h,S}%
\geq2}}t_{h,S}^{2a_{h,S}}\right) \\
&  =2\left(  t_{h,S}{+\frac{t_{h,S}{^{4}}}{{1-t_{h,S}^{2}}}}\right)
\leq2t_{h,S}\left(  1{+\frac{\overline{\tau}{^{3}}}{{1-}\overline{\tau}^{2}}%
}\right)  =3t_{h,S},
\end{align*}
by definition of $\overline{\tau}$. Then:
\begin{gather*}
\sum_{a_{h,S}\in\mathbf{Z}_{p}-\{0\}}\left(  1-2{{\mu(\overline{\mathbf{u}%
}_{h,S})\mu(\overline{\mathbf{v}}_{h,S})+2{\mu(\overline{\mathbf{u}}_{h,S}%
)\mu(\overline{\mathbf{v}}_{h,S})}\cos\left(  {\frac{{2\pi}}{p}}%
a_{h,S}\right)  }}\right)  ^{|M_{h,S}|}\\
\leq3{\exp\left(  {-\frac{{16}}{{lkp^{2}}}{\mu(\overline{\mathbf{u}}_{h,S}%
)\mu(\overline{\mathbf{v}}_{h,S})}(n-(dlk+1))}\right)  }%
\end{gather*}%
\[
\Rightarrow\prod_{h=1}^{|S|}\sum_{a_{h,S}\in\mathbf{Z}_{p}-\{0\}}\left(
1-2{{\mu(\overline{\mathbf{u}}_{h,S})\mu(\overline{\mathbf{v}}_{h,S}%
)+2{\mu(\overline{\mathbf{u}}_{h,S})\mu(\overline{\mathbf{v}}_{h,S})}%
\cos\left(  {\frac{{2\pi}}{p}}a_{h,S}\right)  }}\right)  ^{|M_{h,S}|}%
\]%
\[
\leq{3}^{|S|}\exp\left(  {-\frac{{2\alpha_{S}(n-(dlk+1))}}{{p^{2}}}}\right)
,\quad\text{where }\alpha_{S}=\frac{8}{lk}\overset{|S|}{\sum\limits_{h=1}%
}{{\mu(\overline{\mathbf{u}}_{h,S})\mu(\overline{\mathbf{v}}_{h,S})}}.
\]
Moreover, from (\ref{form5}) and (\ref{zeta}):
\begin{align*}
\sum_{Z^{c}}|\widehat{{P}}_{n}(\mathbf{\alpha})|^{2}  &  \leq\exp\left(
{-\frac{{2\alpha(n-(dlk+1))}}{{p^{2}}}}\right)  \overset{k}{\sum_{|S|=1}%
}\left(
\begin{array}
[c]{c}%
k\\
|S|
\end{array}
\right)  3^{|S|}\\
&  =\left(  4^{k}-1\right)  \exp\left(  {-\frac{{2\alpha(n-(dlk+1))}}{{p^{2}}%
}}\right)  ,
\end{align*}
where $\alpha=\underset{\emptyset\subset S\subset\left\{  {1,...,k}\right\}
}{\min}\alpha_{S}=\frac{8}{lk}\underset{h,S}{\min}{{\mu(\overline{\mathbf{u}%
}_{h,S})\mu(\overline{\mathbf{v}}_{h,S})=-}}\frac{\ln\overline{\tau}%
}{\overline{c}}$. Then, from (\ref{form1}), (\ref{form4}) and (\ref{eps1}):
\[
\Vert P_{n}-U\Vert^{2}\leq\frac{1}{4}\left(  {{{}\varepsilon_{1}(p,n)+}%
}\left(  4^{k}-1\right)  \exp\left(  {-\frac{{2\alpha(n-(dlk+1))}}{{p^{2}}}%
}\right)  \right)  .
\]
Observe that $\displaystyle\underset{p\rightarrow\infty}{\lim}\frac
{{2\alpha(dlk+1)}}{{p^{2}}}=0$ and, if $n\geq cp^{2}$, $\displaystyle
\underset{p\rightarrow\infty}{\lim}\frac{{\varepsilon_{1}(p,n)}}{{\exp}\left(
-2\alpha n/p^{2}\right)  }=0$. Then, if $p$ is sufficiently large:
\[
\Vert P_{n}-U\Vert^{2}\leq\frac{1}{4}\left(  {4^{k}\exp\left(  {-\frac
{{2\alpha n}}{{p^{2}}}}\right)  }\right)  ,
\]
from which
\[
\Vert P_{n}-U\Vert\leq2^{k-1}\exp\left(  {-\frac{{\alpha n}}{{p^{2}}}}\right)
.\ \Box
\]
\medskip

\textbf{Theorem 3.10.}\ Suppose that the matrix $A$ has an eigenvalue
$\lambda\in\mathbf{C},\ |\lambda|>1$ (hence, the matrix ${}^{t}A$ too), that
the support of $\mu$ is not parallel to a proper subspace of $\mathbf{Q}^{k}$
invariant under $A$ and that $\left\|  \mathbf{B}_{n}\right\|  _{\infty}\in
L^{2}$, for all $n\in\mathbf{N}$. Then, there exist $c\in\mathbf{R}^{+}$ and
$N\in\mathbf{N}$ such that, for all $p\in\mathbf{N}$ such that $p>N,\ \gcd
(\det(A),p)=\gcd(\det(B),p)=1$, and for all $n\leq c\ln p$, we have:
\[
\Vert P_{n}-U\Vert\geq{\frac{1}{2}}\eta(p),\qquad\text{where }\lim
_{p\rightarrow\infty}\eta(p)=1.
\]
Consequently, $O(\ln p)$ steps are needed to reach the uniform distribution.

\textbf{Proof.}\ If ${}^{t}A$ has an eigenvalue $\lambda\in\mathbf{C}%
,\ |\lambda|>1$, then:%
\[
\left\|  ^{t}A\right\|  _{\infty}\equiv\underset{\mathbf{x}\in\mathbf{C}%
^{k}-\{\mathbf{0}\}}{\sup}\frac{\left\|  ^{t}A\mathbf{x}\right\|  _{\infty}%
}{\left\|  \mathbf{x}\right\|  _{\infty}}\geq|\lambda|>1.
\]
Let $\mathbf{\alpha}\in\mathbf{Z}_{p}^{k}-\{\mathbf{0}\}$; from the lemmas 2.2
and 2.3:
\begin{gather*}
\Vert P_{n}-U\Vert\geq{\frac{1}{2}}\left|  \widehat{P}_{n}(\mathbf{\alpha
})\right| \\
=\frac{1}{2}\prod\limits_{j=0}^{n-1}\left(  {\sum\limits_{\mathbf{h,i}%
\in\mathbf{Z}^{K}}{\mu(\mathbf{h})\mu(\mathbf{i})\cos\left(  {\frac{{2\pi}}%
{p}\left\langle {\mathbf{h-i},^{t}A^{j}}\mathbf{\alpha}\right\rangle }\right)
}}\right)  ^{1/2}.
\end{gather*}
Since $\displaystyle\cos x\geq1-\frac{{x^{2}}}{2}\,\forall~x\in\mathbf{R}$, we
have:
\begin{equation}
||P_{n}-U||\geq{\frac{1}{2}}\prod\limits_{j=0}^{n-1}{}\left(  {1-\frac{{\rho
}\left\|  ^{t}A\right\|  _{\infty}^{2j}}{{p^{2}}}}\right)  ^{1/2}%
,\quad\text{where }\rho=2\pi^{2}k^{2}\left\|  \mathbf{\alpha}\right\|
_{\infty}^{2}{\sum\limits_{\mathbf{h,i}\in\mathbf{Z}^{K}}{\mu(\mathbf{h}%
)\mu(\mathbf{i})}}\left\|  {{{{\mathbf{h-i}}}}}\right\|  _{\infty}^{2}%
\in\mathbf{R}^{+}. \label{inv}%
\end{equation}
Moreover, $\exists\ d\in\mathbf{R}^{+}$ such that $1-x\geq\exp(-2x)$%
,$\ \forall\ x\in\lbrack0,d]$.

Suppose $n\leq c\ln p$, where $\displaystyle c<{\frac{1}{\ln\left\|
^{t}A\right\|  _{\infty}}}$. Then, $\forall~j=0,1,...,n-1$:
\[
\left\|  ^{t}A\right\|  _{\infty}^{j}<\left\|  ^{t}A\right\|  _{\infty}%
^{n}\leq\left\|  ^{t}A\right\|  _{\infty}^{c\ln p}=p^{c\ln\left\|
^{t}A\right\|  _{\infty}}\Rightarrow\frac{\left\|  ^{t}A\right\|  _{\infty
}^{2j}}{p^{2}}<p^{2\left(  c\ln\left\|  ^{t}A\right\|  _{\infty}-1\right)  }.
\]
Since $\displaystyle\underset{p\rightarrow\infty}{\lim}p^{2\left(
c\ln\left\|  ^{t}A\right\|  _{\infty}-1\right)  }=0$, for sufficiently large
$p$ we can suppose $\displaystyle{\frac{{\rho}\left\|  ^{t}A\right\|
_{\infty}^{2j}}{{p^{2}}}}\in\lbrack0,d]$; hence:
\[
||P_{n}-U||\geq\frac{1}{2}\prod\limits_{j=0}^{n-1}\exp\left(  {-\frac{{\rho
}\left\|  ^{t}A\right\|  _{\infty}^{2j}}{{p^{2}}}}\right)  =\frac{1}{2}%
\exp\left(  {-\frac{\rho}{{p^{2}}}\sum\limits_{j=0}^{n-1}{\left(  \left\|
^{t}A\right\|  _{\infty}^{2}\right)  ^{j}}}\right)
\]%
\[
=\frac{1}{2}\exp\left(  {-\frac{\rho}{\left\|  ^{t}A\right\|  _{\infty}%
^{2}{-1}}\cdot\frac{\left\|  ^{t}A\right\|  _{\infty}^{2n}{-1}}{{p^{2}}}%
}\right)
\]%
\[
\geq\frac{1}{2}\exp\left(  {-\frac{\rho}{\left\|  ^{t}A\right\|  _{\infty}%
^{2}{-1}}\cdot p^{2\left(  {c\ln\left\|  ^{t}A\right\|  _{\infty}-1}\right)
}+\frac{\rho}{{\left(  \left\|  ^{t}A\right\|  _{\infty}^{2}{-1}\right)
p^{2}}}}\right)  \equiv\frac{1}{2}\eta(p).
\]

By definition of $c$, we have the thesis.\ $\Box$\medskip

\textbf{Theorem 3.11.}\ Suppose that the matrix $A$ has an eigenvalue
$\lambda\in\mathbf{C}$ such that $\lambda^{l}=1$, for some $l\in
\mathbf{N}^{\ast}$ (hence, the matrix $^{t}A$ too), that the support of $\mu$
is not parallel to a proper subspace of $\mathbf{Q}^{k}$ invariant under $A$
and that $\left\|  \mathbf{B}_{n}\right\|  _{\infty}\in L^{2}$, for all
$n\in\mathbf{N}$. Then, there exist $\gamma\in\mathbf{R}^{+}$ and
$N\in\mathbf{N}$ such that, for all $p\in\mathbf{N}$ such that $p>N,\ \gcd
(\det(A),p)=\gcd(\det(B,p)=1$, we have:
\[
||P_{n}-U||\geq\frac{1}{2}\exp\left(  {-\frac{{\gamma n}}{{p^{2}}}}\right)  .
\]
Consequently, $O\left(  p^{2}\right)  $ steps are needed to reach the uniform distribution.

\textbf{Proof.}\ The assumption on $\lambda$ implies ${^{t}A^{l}}%
\mathbf{x}=\mathbf{x}$, for some $\mathbf{x}\in\mathbf{C}^{k}-\{\mathbf{0\}}$,
and so $({^{t}A^{l}-I)}\mathbf{x}=\mathbf{0}$, which implies $\mathbf{x}%
\in\mathbf{Q}^{k}-\{\mathbf{0\}}$; then $\exists$ $\mathbf{\alpha}%
\in\mathbf{Z}^{k}-\{\mathbf{0\}}$ such that ${^{t}A^{l}}\mathbf{\alpha
}=\mathbf{\alpha}$. $\forall$ $p>\left\|  \mathbf{\alpha}\right\|  _{\infty}$,
we can suppose $\mathbf{\alpha}\in\mathbf{Z}_{p}^{k}-\{\mathbf{0\}}$,
${^{t}A^{l}}\mathbf{\alpha}=\mathbf{\alpha}$ (mod $p$); then, $\forall
~j\in\mathbf{N}$, $\exists~i\in\{0,1,...,l-1\}$ such that ${^{t}A^{j}%
}\mathbf{\alpha}={^{t}A^{i}}\mathbf{\alpha}$ (mod $p$).

By proceeding as in the proof of the previous theorem, we obtain the following
formula, analogous to (\ref{inv}):
\[
||P_{n}-U||\geq\frac{1}{2}\left(  {1-\frac{\gamma}{{p^{2}}}}\right)  ^{n/2},
\]
where $\displaystyle\gamma=2\pi^{2}k^{2}\left\|  \mathbf{\alpha}\right\|
_{\infty}^{2}\underset{i\in\{0,1,...,l-1\}}{\max}\left\|  ^{t}A\right\|
_{\infty}^{2i}{\sum\limits_{\mathbf{h,i}\in\mathbf{Z}^{K}}{\mu(\mathbf{h}%
)\mu(\mathbf{i})}}\left\|  {{{{\mathbf{h-i}}}}}\right\|  _{\infty}^{2}%
\in\mathbf{R}^{+}$. Finally, $\forall~p$ sufficiently large:%
\[
||P_{n}-U||\geq\frac{1}{2}\exp\left(  {-\frac{{\gamma n}}{{p^{2}}}}\right)
.\ \Box
\]

\section{Problems for further study}

A natural problem to study is the generalization of the recursion (\ref{iter})
to the analogous recursion \ in $\mathbf{R}^{k}$ reduced modulo $c$, for some
real number $c$. In this case, the idea is to use the Fourier transform on
$\mathbf{R}^{k}$ instead of in $\mathbf{Z}^{k}$ and then generalize the lemmas
in Section 2. The expectation is to prove the convergence in law of the Markov
chain and to estimate the rate of convergence to the uniform distribution on
some subset of $\mathbf{R}^{k}$: the set where the chain ranges. However, this
set can be different from $\mathbf{R}^{k}$ (mod $c$) and it can be also
countable. In order to establish it and to develop the theory, some changes of
the results in Section 2 are needed: for example, Lemma 2.1 (upper bound
lemma) is not valid in the continuous context and it must be modified.
Conversely, results as Lemma 2.5 seem useful also in the study of the high
powers of the real matrix $A$ in the modified recursion (\ref{iter}).\

\medskip

Claudio Asci, Dipartimento di Matematica e Informatica, {Universit\`{a} degli
Studi di Trieste, }Via Valerio 12/1, {34127 Trieste, Italy}

E-mail: asci@dmi.units.it
\end{document}